# CONTACT-BOUNDARY VALUE PROBLEM IN THE NON-CLASSICAL TREATMENT FOR ONE PSEUDOPARABOLIC EQUATION

**I.G.Mamedov**

*A.I.Huseynov Institute of Cybernetics of NAS of Azerbaijan. Az 1141, Azerbaijan, Baku st. B. Vahabzade, 9*
*E-mail: ilgar-mammadov@rambler.ru*

## Abstract

*In the paper, the contact - boundary value problem with non-classical conditions not requiring agreement conditions is considered for a pseudoparabolic equation. The equivalence of these conditions is substantiated in the case if the solution of the solution of the stated problem is sought in S.L.Sobolev isotropic space.*

## Problem statement

Consider equation

$$(V_{4,4}u)(x) \equiv \sum_{i_1=0}^{4}\sum_{i_2=0}^{4} a_{i_1,i_2}(x) D_1^{i_1} D_2^{i_2} u(x) = Z_{4,4}(x) \in L_p(G), \qquad (1)$$

where $a_{4,4}(x) \equiv 1$.

Here $u(x) \equiv u(x_1, x_2)$ is a desired function determined on $G$; $a_{i_1,i_2}(x)$ are the given measurable functions on $G = G_1 \times G_2$, where $G_k = (0, h_k)$, $k = 1,2$. $Z_{4,4}(x)$ is a given measurable function on $G$; $D_k = \partial / \partial x_k$ is a generalized differentiation operator in S.L.Sobolev sense, $D_k^0$ is an identity transformation operator.

Equation (1) is a hyperbolic equation that has two real characteristics $x_1 = const, x_2 = const$, the first and second one of which is four-fold. Therefore, we can consider equation (1) in some sense as a pseudoparabolic equation [1]. This equation is a generalization of the equation of thin spherical shell bending [2, p.258].

In the paper, we consider equation (1) in the general case when the coefficients $a_{i_1,i_2}(x)$ are non smooth functions satisfying the following conditions:

$$a_{i_1,i_2}(x) \in L_p(G), \; i_1 = \overline{0,3} \; i_2 = \overline{0,3};$$

$$a_{4,i_2}(x) \in L_{\infty,p}^{x_1,x_2}(G), \; i_2 = \overline{0,3};$$

$$a_{i_1,4}(x) \in L_{p,\infty}^{x_1,x_2}(G), \; i_1 = \overline{0,\tilde{3}}$$

Therewith, the important principal moment is that the equation under consideration has discontinuous coefficients that satisfy only some *p*-integrability and boundedness conditions, i.e. the considered pseudoparabolic differential operator has no traditional conjugation operator.

Under these conditions, we'll look for the solution $u(x)$ of equation (1) in S.L.Sobolev isotropic space

$$W_p^{(4,4)}(G) \equiv \{u(x): D_1^{i_1} D_2^{i_2} u(x) \in L_p(G), i_1 = \overline{0,4}, i_2 = \overline{0,4}\}$$

where $1 \leq p \leq \infty$.

For equation (1) we can give the contact-boundary conditions of the classic form as follows: (see fig.1 ):

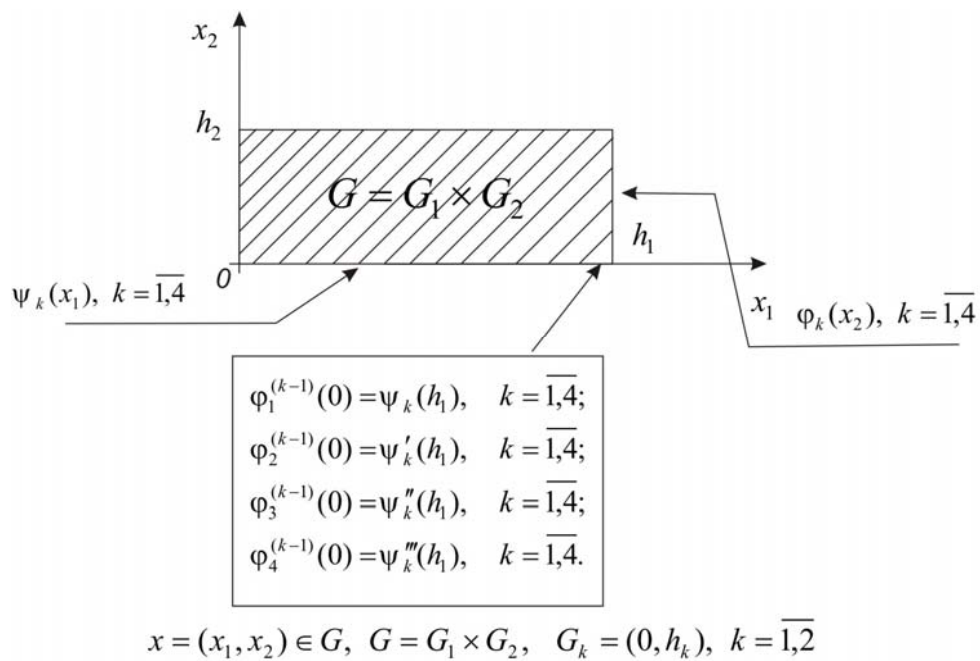

Fig. 1. Geometric interpretation of classic contact-boundary conditions.

$$\begin{cases}
u(x_1,x_2)\big|_{x_1=h_1}=\varphi_1(x_2); & u(x_1,x_2)\big|_{x_2=0}=\psi_1(x_1); \\
\dfrac{\partial u(x_1,x_2)}{\partial x_1}\bigg|_{x_1=h_1}=\varphi_2(x_2); & \dfrac{\partial u(x_1,x_2)}{\partial x_2}\bigg|_{x_2=0}=\psi_2(x_1); \\
\dfrac{\partial^2 u(x_1,x_2)}{\partial x_1^2}\bigg|_{x_1=h_1}=\varphi_3(x_2); & \dfrac{\partial^2 u(x_1,x_2)}{\partial x_2^2}\bigg|_{x_2=0}=\psi_3(x_1); \\
\dfrac{\partial^3 u(x_1,x_2)}{\partial x_1^3}\bigg|_{x_1=h_1}=\varphi_4(x_2); & \dfrac{\partial^3 u(x_1,x_2)}{\partial x_2^3}\bigg|_{x_2=0}=\psi_4(x_1);
\end{cases} \quad (2)$$

where $\varphi_k(x_2),\ \psi_k(x_1),\ k=\overline{1,4}$ are the given measurable functions on $G$. It is obvious that in the case of conditions (2), in addition to the conditions

$$\varphi_k(x_2)\in W_p^{(4)}(G_2),\ \psi_k(x_1)\in W_p^{(4)}(G_1)$$

the given functions should satisfy also the following agreement conditions:

$$\begin{cases}
\varphi_1(0)=\psi_1(h_1); & \varphi_2(0)=\psi_1'(h_1); \\
\varphi_1'(0)=\psi_2(h_1); & \varphi_2'(0)=\psi_2'(h_1); \\
\varphi_1''(0)=\psi_3(h_1); & \varphi_2''(0)=\psi_3'(h_1); \\
\varphi_1'''(0)=\psi_4(h_1); & \varphi_2'''(0)=\psi_4'(h_1); \\
\varphi_3(0)=\psi_1''(h_1); & \varphi_4(0)=\psi_1'''(h_1); \\
\varphi_3'(0)=\psi_2''(h_1); & \varphi_4'(0)=\psi_2'''(h_1); \\
\varphi_3''(0)=\psi_3''(h_1); & \varphi_4''(0)=\psi_3'''(h_1); \\
\varphi_3'''(0)=\psi_4''(h_1); & \varphi_4'''(0)=\psi_4'''(h_1).
\end{cases} \quad (3)$$

Obviously, conditions (2) are close to the boundary conditions of the Goursat problem from [3-6].

Consider the following non-classical boundary conditions:

$$V_{i_1,i_2}u \equiv D_1^{i_1}D_2^{i_2}u(h_1,0)=Z_{i_1,i_2}\in R,\ i_k=\overline{0,3},\ k=\overline{1,2};$$
$$(V_{4,i_2}u)(x_1)\equiv D_1^4 D_2^{i_2}u(x_1,0)=Z_{4,i_2}(x_1)\in L_p(G_1),\ i_2=\overline{0,3}; \quad (4)$$
$$(V_{i_1,4}u)(x_2)\equiv D_1^{i_1}D_2^4 u(h_1,x_2)=Z_{i_1,4}(x_2)\in L_p(G_2),\ i_1=\overline{0,3}.$$

If the function $u \in W_p^{(4,4)}(G)$ is a solution of the classical form contact boundary value problem (1), (2), then it is also a solution of problem (1), (4) for $Z_{i_1,i_2}$, determined by the following equalities:

$$Z_{0,0} = \varphi_1(0) = \psi_1(h_1); \quad Z_{0,1} = \varphi_1'(0) = \psi_2(h_1); \quad Z_{1,0} = \varphi_2(0) = \psi_1'(h_1);$$

$$Z_{1,1} = \varphi_2'(0) = \psi_2'(h_1); \quad Z_{2,0} = \varphi_3(0) = \psi_1''(h_1); \quad Z_{2,1} = \varphi_3'(0) = \psi_2''(h_1);$$

$$Z_{3,0} = \varphi_4(0) = \psi_1'''(h_1); \quad Z_{3,1} = \varphi_4'(0) = \psi_2'''(h_1); \quad Z_{0,2} = \varphi_1''(0) = \psi_3(h_1);$$

$$Z_{0,3} = \varphi_1'''(0) = \psi_4(h_1); \quad Z_{1,2} = \varphi_2''(0) = \psi_3'(h_1); \quad Z_{1,3} = \varphi_2'''(0) = \psi_4'(h_1);$$

$$Z_{2,2} = \varphi_3''(0) = \psi_3''(h_1); \quad Z_{2,3} = \varphi_3'''(0) = \psi_4''(h_1); \quad Z_{3,2} = \varphi_4''(0) = \psi_3'''(h_1);$$

$$Z_{3,3} = \varphi_4'''(0) = \psi_4'''(h_1); \quad Z_{4,0}(x_1) = \psi_1^{(IV)}(x_1); \quad Z_{4,2}(x_1) = \psi_3^{(IV)}(x_1);$$

$$Z_{4,1}(x_1) = \psi_2^{(IV)}(x_1); \quad Z_{4,3}(x_1) = \psi_4^{(IV)}(x_1); \quad Z_{0,4}(x_2) = \varphi_1^{(IV)}(x_2);$$

$$Z_{2,4}(x_2) = \varphi_3^{(IV)}(x_2); \quad Z_{1,4}(x_2) = \varphi_2^{(IV)}(x_2); \quad Z_{3,4}(x_2) = \varphi_4^{(IV)}(x_2)$$

It is easily proved that the inverse is also true. In other words, if the function $u \in W_p^{(4,4)}(G)$ is a solution of problem (1), (4), it is also a solution of problem (1), (2) for the following functions:

$$\varphi_1(x_2) = Z_{0,0} + x_2 Z_{0,1} + \frac{x_2^2}{2!} Z_{0,2} + \frac{x_2^3}{3!} Z_{0,3} + \int_0^{x_2} \frac{(x_2-\tau)^3}{3!} Z_{0,4}(\tau) d\tau; \quad (5)$$

$$\varphi_2(x_2) = Z_{1,0} + x_2 Z_{1,1} + \frac{x_2^2}{2!} Z_{1,2} + \frac{x_2^3}{3!} Z_{1,3} + \int_0^{x_2} \frac{(x_2-\xi)^3}{3!} Z_{1,4}(\xi) d\xi; \quad (6)$$

$$\varphi_3(x_2) = Z_{2,0} + x_2 Z_{2,1} + \frac{x_2^2}{2!} Z_{2,2} + \frac{x_2^3}{3!} Z_{2,3} + \int_0^{x_2} \frac{(x_2-\eta)^3}{3!} Z_{2,4}(\eta) d\eta; \quad (7)$$

$$\varphi_4(x_2) = Z_{3,0} + x_2 Z_{3,1} + \frac{x_2^2}{2!} Z_{3,2} + \frac{x_2^3}{3!} Z_{3,3} + \int_0^{x_2} \frac{(x_2-\nu)^3}{3!} Z_{3,4}(\nu) d\nu; \quad (8)$$

$$\psi_1(x_1) = Z_{0,0} + (x_1-h_1) Z_{1,0} + \frac{(x_1-h_1)^2}{2!} Z_{2,0} +$$

$$+ \frac{(x_1-h_1)^3}{3!} Z_{3,0} + \int_{h_1}^{x_1} \frac{(x_1-\lambda)^3}{3!} Z_{4,0}(\lambda) d\lambda; \quad (9)$$

$$\psi_2(x_1) = Z_{0,1} + (x_1 - h_1)Z_{1,1} + \frac{(x_1 - h_1)^2}{2!}Z_{2,1} +$$
$$+ \frac{(x_1 - h_1)^3}{3!}Z_{3,1} + \int_{h_1}^{x_1} \frac{(x_1 - \mu)^3}{3!} Z_{4,1}(\mu)d\mu; \qquad (10)$$

$$\psi_3(x_1) = Z_{0,2} + (x_1 - h_1)Z_{1,2} + \frac{(x_1 - h_1)^2}{2!}Z_{2,2} +$$
$$+ \frac{(x_1 - h_1)^3}{3!}Z_{3,2} + \int_{h_1}^{x_1} \frac{(x_1 - \rho)^3}{3!} Z_{4,2}(\rho)d\rho; \qquad (11)$$

$$\psi_4(x_1) = Z_{0,3} + (x_1 - h_1)Z_{1,3} + \frac{(x_1 - h_1)^2}{2!}Z_{2,3} +$$
$$+ \frac{(x_1 - h_1)^3}{3!}Z_{3,3} + \int_{h_1}^{x_1} \frac{(x_1 - \sigma)^3}{3!} Z_{4,3}(\sigma)d\sigma; \qquad (12)$$

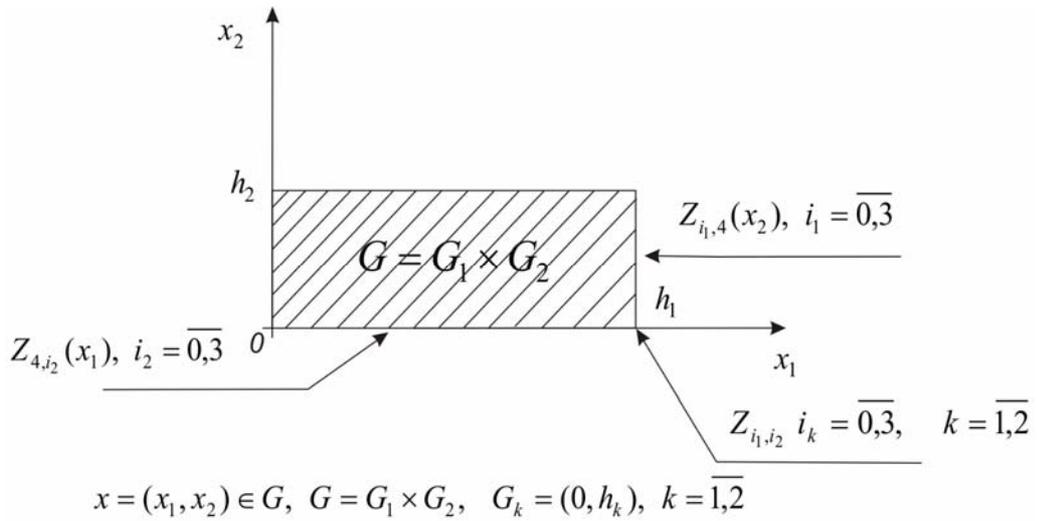

**Fig. 2. Geometric interpretation of contact-boundary value conditions in non-classical treatment.**

Note that the functions (5)-(12) possess an important property, more exactly, the agreement conditions for all $Z_{i_1,i_2}$, possessing the above-mentioned properties are fulfilled for them. Therefore, equalities (5)-(12) may be considered as a general form of all the functions $\varphi_k(x_2), \psi_k(x_1)$, $k=\overline{1,4}$ satisfying agreement conditions (3).

So, the classic form contact-boundary value problem (1), (2) and in non-classic treatment (1), (4) (see. Fig.2) are equivalent in the general case. However, the contact-boundary value problem in non-classic treatment (1), (4) is more natural by the statement than problem (1), (2). This is connected with the fact that in the statement of problem (1), (4), the right sides of boundary conditions don't require additional conditions of agreement type. Note that some boundary value problems in non-classic treatments were considered in the author's papers [7-9].